\documentclass{amsart}

\usepackage{graphicx}

\theoremstyle{definition}

\theoremstyle{remark}

\numberwithin{equation}{section}

\begin{document}

\title{The \(\zeta\)-regularized product over all primes}
\author{V.V. Smirnov} 
\email{ixorg@ya.ru} 

\begin{abstract}
In this paper we prove that the \(\zeta\)-regularized product over all primes is \(\pi e^{\mu}\), where \(\mu\) is
closely related  with the non-trivial zeros of the \(\zeta(s)\).
\end{abstract}

\maketitle

\section{Introduction.}
The Chebyshev theta function \(\theta(x)\) is defined by
\begin{equation} 
\theta(x)=\sum_{p\le x}\log p,
\end{equation}
where \(p\) runs over primes. There is a simple relationship between  \(\theta(x)\) and the prime-counting function \(\pi(x)\):
\begin{equation} 
\pi(x)=\frac{\theta(x)}{\log x}+\int_2^x \frac{\theta(t)}{t\log^2 t}dt.  
\end{equation}
Indeed, let \(p_k\in\mathbb{P}\) and let \(a_k=\sum_{j=1}^k\log p_j\) so \(a_k-a_{k-1}=\log⁡   p_k\), hence we have
\begin{eqnarray} 
\begin{split} 
\int_{2}^{x} \frac{\theta(t)}{t\log^2 t}dt&=\sum_{j=1}^{k-1} \int_{p_j}^{p_{j+1}} \frac{a_j}{t\log^2 t}dt+\int_{p_k}^{x} \frac{a_k}{t\log^2 t}dt
\nonumber\\
&=\sum_{j=1}^{k-1} a_j\left(\frac{1}{\log p_j}-\frac{1}{\log p_{j+1}}\right)+\frac{a_k}{\log p_k}-\frac{a_k}{\log x}
\nonumber\\
&=\frac{a_1}{\log p_1}+\sum_{j=2}^{k} \frac{a_j-a_{j-1}}{\log p_j}-\frac{a_k}{\log x}\nonumber\\
&=1+\sum_{j=2}^{k}1-\frac{\theta(x)}{\log x}
=\pi(x)-\frac{\theta(x)}{\log x}.
\end{split} 
\end{eqnarray}
The second Chebyshev function \(\psi(x)\) is defined similarly, with the sum extending over all prime powers not exceeding \(x\):
\begin{equation} 
\psi(x)=\sum_{p^k\le x}\log p=\sum_{n\le x}\Lambda(n),
\end{equation}
where \(\Lambda(n)\) is the von Mangoldt function.  \(\psi(x)\) can be related to the prime-counting function as follows
\begin{equation} 
\Pi(x)=\sum_{n\le x}\Lambda(n)\int_{n}^{x} \frac{dt}{t\log^2 t}+\frac{1}{\log x}\sum_{n\le x}\Lambda(n)=\frac{\psi(x)}{\log x}+\int_{2}^{x} \frac{\psi(t)}{t\log^2 t}dt.
\end{equation}
\begin{equation} 
\Pi(x)=\pi(x)+\frac{1}{2}\pi(x^{1/2})+\frac{1}{3}\pi(x^{1/3})+\cdots.
\end{equation}
Also, \(\psi(x)\) has an explicit expression as a sum over the nontrivial zeros of the Riemann zeta function [1, p.104]:
\begin{equation} 
\psi_0(x)=\frac{1}{2}\left(\sum_{n\le x}\Lambda(n)+\sum_{n<x}\Lambda(n)\right)=x-\sum_\rho \frac{x^\rho}{\rho}-\frac{\zeta^\prime(0)}{\zeta(0)}-\sum_{k=1}^{\infty}\frac{x^{-2k}}{-2k},
\end{equation}
where \(\rho\) runs over the nontrivial zeros of the zeta function. The terms on the right side are arranged by decreasing magnitude. Thus
\begin{equation} 
\sum_{p^k\le x}\log p \sim x
\end{equation}
and the location of the nontrivial zeros of \(\zeta(s)\) determines the largest error term. Since the zeros are symmetric about the line  \(\Re e(s)=\frac{1}{2} \),
 and since \(|x^\rho |=x^{\Re e(\rho)}\), the error term is as small as possible if and only if \(\Re e(\rho)=\frac{1}{2}\) for all  \(\rho\).

\section{the constant.}
The Riemann \(\xi\) function is defined as:
\begin{equation} 
 \xi (s)=\Pi \left(\frac{s}{2}\right)(s-1)\pi^{-s/2}\zeta(s)
\end{equation}
for  \(s\in\mathbb{C}\),  where 
\begin{equation} 
 \Pi(s)=s\Gamma(s)=\Gamma(s+1)=\int_0^\infty e^{-t}t^s dt.
\end{equation}
\(\xi(s)\) is an entire function -- that is, an analytic function of \(s\) which is defined for all values of \(s\)
 -- and the functional equation of the zeta function is equivalent to  \(\xi(s)=\xi(1-s)\).
 Xi-function can also be expanded as infinite product
\begin{equation} 
 \xi(s)=\xi(0)\prod_\rho \left(1-\frac{s}{\rho}\right),
\end{equation}
where \(\rho\) ranges over the roots of the equation \(\xi(\rho)=0\). Since the logarithmic derivative of \(\xi(s)\) is on the one hand
\begin{equation} 
 \sum_\rho \frac{d}{ds}\log \left(1-\frac{s}{\rho}\right)=\sum_\rho \frac{1}{s-\rho}
\end{equation}
and on the other hand
\begin{equation} 
 \frac{d}{ds}\log\Pi\left(\frac{s}{2}\right)-\frac{1}{2}\log\pi+\frac{1}{s-1}+\frac{\zeta^\prime (s)}{\zeta(s)}
\end{equation}
the sum of the series \(\sum 1/\rho\) is the value of
\begin{equation} 
 -\frac{1}{2} \frac{\Pi^\prime (s/2)}{\Pi(s/2)} +\frac{1}{2}\log\pi+\frac{1}{1-s}-\frac{\zeta^\prime (s)}{\zeta(s)}
\end{equation}
at \(s=0\). Now logarithmic differentiation of the product formula for \(\Pi(s)\)
\begin{equation} 
 \Pi(s)=\prod_{n=1}^\infty \frac{n^{1-s}(n+1)^s}{s+n}=\prod_{n=1}^\infty \left(1+\frac{s}{n}\right)^{-1} \left(1+\frac{1}{n}\right)^s
\end{equation}
gives
\begin{equation} 
 \frac{\Pi^\prime(s)}{\Pi(s)}=\sum_{n=1}^\infty \left[-\frac{1}{s+n}-\log n + \log(n+1)\right],
\end{equation}
\begin{equation} 
 -\frac{\Pi^\prime(0)}{\Pi(0)}=\lim\limits_{n\to\infty} \left[1+\frac{1}{2}+\frac{1}{3}+\cdots+\frac{1}{n}-\log(n+1)\right].
\end{equation}
The number on the right side of this equation is by definition \textit{Euler–-Mascheroni constant} and is traditionally denoted \(\gamma\). Thus [3, \textsection 41],[5],[6]
\begin{equation} 
 \eta:=\sum_\rho \frac{1}{\rho(1-\rho)}=2\sum_\rho \frac{1}{\rho}=\gamma+2-\log 4\pi=0.04619\ldots
\end{equation}

\section{The \(\zeta\)-regularized product over all primes.}
We define the \textit{prime zeta function}  \(P(s)\), \(s=\sigma+i\tau\), through
\begin{equation} 
 P(s)=\sum_p p^{-s}
\end{equation}
with the summation performed over all primes  \(p\). The series converges absolutely when \(\sigma>1\).
But we also know that [4, p. 70]
\begin{equation} 
 e^{p^{-s}}=\prod_{n=1}^{+\infty} (1-p^{-ns})^{-\frac{\mu(n)}{n}},
\end{equation}
where \(\mu(n)\) is the M\"obius function. So it follows that
\begin{equation} 
 e^{P(s)}=\prod_p e^{p^{-s}}=\prod_p \prod_{n=1}^{+\infty} (1-p^{-ns})^{-\frac{\mu(n)}{n}}=\prod_{n=1}^{+\infty} \zeta(ns)^{\frac{\mu(n)}{n}}.
\end{equation}
Taking the logarithmic derivative we get
\begin{equation} 
 P^\prime(s)=\sum_{n=1}^{+\infty} \frac{\mu(n)}{n} \frac{n\zeta^\prime(ns)}{\zeta(ns)}=\sum_{n=1}^{+\infty} \mu(n) \frac{\zeta^\prime(ns)}{\zeta(ns)}.
\end{equation}
Now consider the integral
\begin{eqnarray} 
\begin{split} 
\int_0^\infty \frac{xt^{x-1}}{e^t+1} dt &=\int_0^\infty \frac{xt^{x-1}}{1+e^{-t}}e^{-t} dt\\
&=x\sum_{k=1}^\infty(-1)^{k-1}\int_0^\infty t^{x-1}e^{-kt} dt\\
&=x\sum_{k=1}^\infty(-1)^{k-1}k^{-x}\int_0^\infty t^{x-1}e^{-t} dt\\
&=(1-2^{1-x})\zeta(x)\Gamma(x+1).
\end{split}
\end{eqnarray}
Integration by parts gives
\begin{equation} 
 \lim_{x\to0^+}\int_0^\infty\frac{xt^{x-1}}{e^t+1} dt = \lim_{x\to0^+}\int_0^\infty\frac{t^xe^t}{(e^t+1)^2} dt =\int_1^\infty\frac{du}{(u+1)^2}=\frac{1}{2}.
\end{equation}
Making \(x\to0\) in (3.5) and combining with (3.6), we get \(\zeta(0)=-\frac{1}{2}\). 
We also know that \(\zeta\)-regularization gives [4, p.71]
\begin{equation} 
\sum_{n=1}^\infty \mu(n)= -2.
\end{equation}
Therefore
\begin{equation} 
 P^\prime(0)=\left(\sum_{n=1}^{+\infty} \mu(n)\right)\frac{\zeta^\prime(0)}{\zeta(0)}=\frac{1}{\zeta(0)} \frac{\zeta^\prime(0)}{\zeta(0)}=-2\log2\pi.
\end{equation}
Finally, using (2.10), (3.1) and (3.8) we obtain (3.9), where the product is taken over all primes  \(p\) and \(\mu=2+\gamma-\eta\).
\begin{equation} 
 \prod p = \pi e^{\mu} 
\end{equation}

\bibliographystyle{amsplain}

\end{document}